\documentclass{article}
\usepackage[utf8]{inputenc}
\usepackage[T1]{fontenc}
\usepackage{lmodern}
\usepackage[english]{babel}
\usepackage{amssymb}
\usepackage[geometry]{ifsym}
\usepackage{graphicx}
\usepackage{subfig}
\usepackage[normalem]{ulem}

\newtheorem{lemma}{Lemma}
\newtheorem{theorem}{Theorem}
\newtheorem{corollary}{Corollary}

\newcommand{\DONOTTEX}[1]{}

\newcommand{\qed}{\hspace*{\fill} $\Box$}

\newcommand{\stareq}{=^\star}
\newcommand{\stargeq}{\geq^\star}
\newcommand{\starleq}{\leq^\star}
\newcommand{\starsubseteq}{\subseteq^\star}

\graphicspath{{./figures/}}

\usepackage[usenames]{color} 
\definecolor{hellblau}{rgb}{0.2,0.4,1} 
\definecolor{dunkelblau}{rgb}{0,0,0.8}
\definecolor{dunkelgruen}{rgb}{0,0.5,0}
\usepackage[
	pdftex,
	colorlinks,
	linkcolor=dunkelblau,
	urlcolor=dunkelblau,
	citecolor=dunkelgruen,
	bookmarks=true,
	linktocpage=true,
	pdfsubject={}
]{hyperref} 
\begin{document}

\title{{\bf More on foxes}}
\date{}
\author{Matthias Kriesell
		\and Jens M. Schmidt
		}

\maketitle

\begin{abstract}
  \setlength{\parindent}{0em}
  \setlength{\parskip}{1.5ex}

  An edge in a $k$-connected graph $G$ is called {\em $k$-contractible} if
  the graph $G/e$ obtained from $G$ by contracting $e$ is $k$-connected.
  Generalizing earlier results on $3$-contractible edges in spanning trees of $3$-connected graphs, we prove that (except for the graphs $K_{k+1}$ if $k \in \{1,2\}$)
  (a) every spanning tree of a $k$-connected triangle free graph has two $k$-contractible edges,
	(b) every spanning tree of a $k$-connected graph of minimum degree at least $\frac{3}{2}k-1$ has two $k$-contractible edges,
  (c) for $k>3$, every DFS tree of a $k$-connected graph of minimum degree at least $\frac{3}{2}k-\frac{3}{2}$ has two $k$-contractible edges,
  (d) every spanning tree of a cubic $3$-connected graph nonisomorphic to $K_4$ has at least $\frac{1}{3}|V(G)|-1$ many $3$-contractible edges, and
  (e) every DFS tree of a $3$-connected graph nonisomorphic to $K_4$, the prism, or the prism plus a single edge has two 3-contractible edges.
  We also discuss in which sense these theorems are best possible.
          

  
\end{abstract}

\maketitle

\section{Introduction}

All graphs throughout are assumed to be finite, simple, and undirected.
For terminology not defined here we refer to~\cite{Diestel2010} or~\cite{BondyMurty2007}.
A graph is called {\em $k$-connected} ($k \geq 1$) if $|V(G)|>k$ and $G-T$ is connected for all $T \subseteq V(G)$ with $|T|<k$.
Let $\kappa(G)$ denote the {\em connectivity} of $G$, that is, the largest $k$ such that $G$ is $k$-connected.
A set $T \subseteq V(G)$ is called a {\em smallest separating set} if $|T|=\kappa(G)$ and $G-T$ is disconnected.
By $\mathfrak{T}(G)$ we denote the set of all smallest separating sets of $G$.
An edge $e$ of a $k$-connected graph $G$ is called {\em $k$-contractible} if the graph $G/e$
obtained from $G$ by {\em contracting} $e$, that is, identifying its endvertices and simplifying the result, is $k$-connected.
No edge in $K_{k+1}$ is $k$-contractible, whereas all edges in $K_\ell$ are if $\ell \geq k+2$,
and it is well-known and straightforward to check that, for a noncomplete $k$-connected graph $G$,
an edge $e$ is not $k$-contractible if and only if $\kappa(G)=k$ and $V(e) \subseteq T$ for some ${\mathfrak T}(G)$.

\begin{figure}[htb]
	\centering
	\subfloat[The prism $K_3 \times K_2$, which is no fox. Dashed edges are 3-contractible. Fat edges depict a DFS tree that contains exactly one 3-contractible edge (namely, $e$).]{
		\includegraphics[scale=0.5]{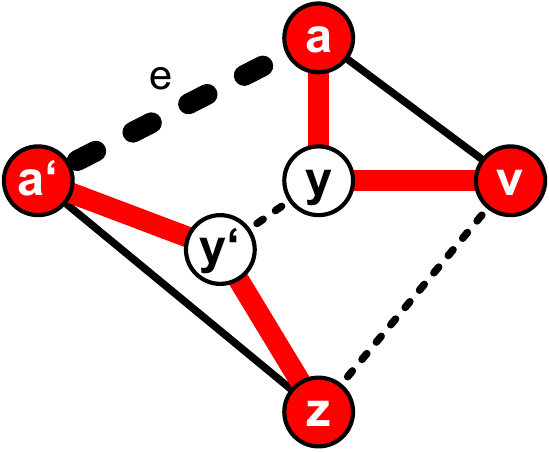}
		\label{fig:6vertices1ContractibleEdgeBW2}
	}
	\hfill
	\subfloat[The fox $(K_3 \times K_2)^+$ and a DFS-tree of it containing exactly one 3-contractible edge.]{
		\includegraphics[scale=0.5]{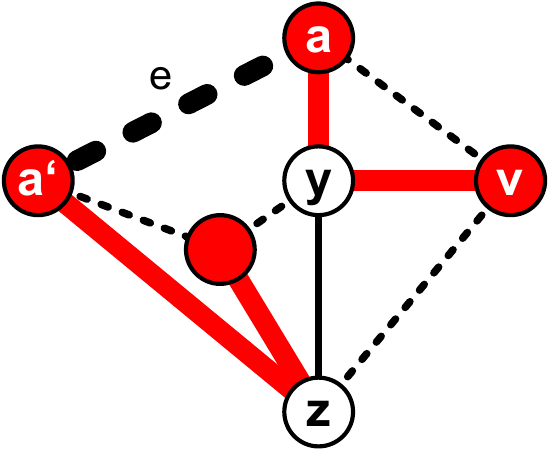}
		\label{fig:6vertices1ContractibleEdgeBW}
	}
	\hfill
	\subfloat[An infinite family of foxes obtained by enlarging the lower left part horizontally. Every fox in this family has a DFS tree (fat edges) that contains exactly two 3-contractible edges.]{
		\includegraphics[scale=0.5]{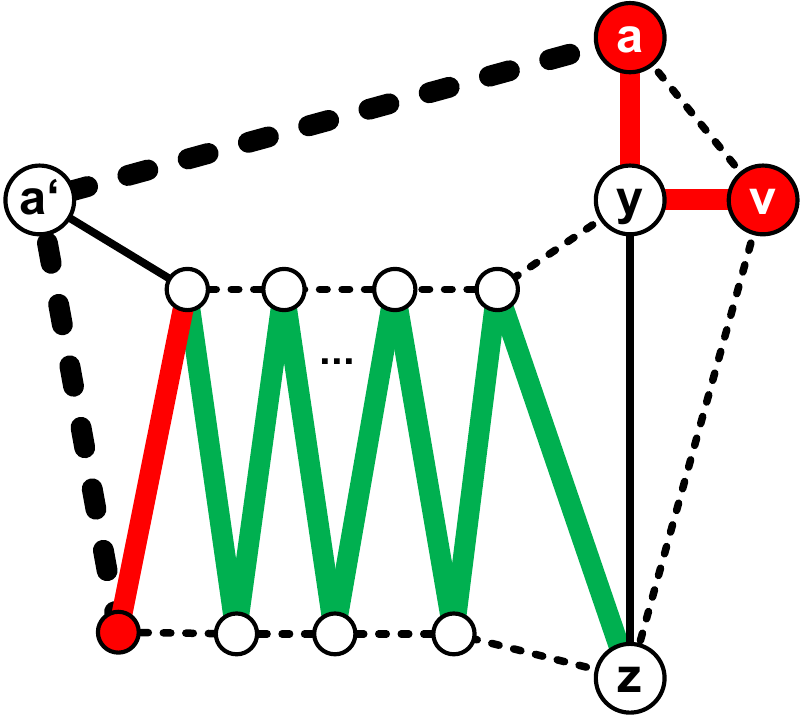}
		\label{fig:LargeDFSTreeWith2ContractibleEdges}
	}
	\caption{The two exceptional graphs of Theorem~\ref{T1} on six vertices and an infinite family of foxes showing that Theorem~\ref{T1} is sharp.}
	\label{fig:Tightness}
\end{figure}

There is a rich literature dealing with the distribution of $k$-contractible edges in $k$-connected graphs (see the surveys~\cite{Kriesell2002,Kriesell2012}),
with a certain emphasis on the case $k=3$. In~\cite{ElmasryMehlhornSchmidt2013}, $3$-connected graphs that admit a spanning tree without any $3$-contractible edge
have been introduced; these were called {\em foxes} (see Figure~\ref{fig:Tightness}). For example, every wheel $G$ is a fox, which is certified by the spanning star $Q$ that is centered at the hub of the wheel.
However, $Q$ is as far from being a \emph{DFS} (depth-first search) tree as it can be, and one could ask if the property of being a fox can be certified by some DFS tree at all.
The answer is no, as it has been shown in~\cite{ElmasryMehlhornSchmidt2013}
that every DFS tree of every $3$-connected graph nonisomorphic to $K_4$ does contain a $3$-contractible edge.
Here we generalize the latter result as follows.

\begin{theorem}
  \label{T1}
  Every DFS tree of every $3$-connected graph nonisomorphic to $K_4$,
  the prism $K_3 \times K_2$, or the unique graph $(K_3 \times K_2)^+$ obtained from $K_3 \times K_2$ by adding a single edge
  contains at least two $3$-contractible edges.
\end{theorem}

Theorem~\ref{T1} is best possible in the sense that there is an infinite class of $3$-connected graphs admitting a DFS tree with only two $3$-contractible edges (see Figure~\ref{fig:LargeDFSTreeWith2ContractibleEdges}).

Our proofs are based on methods introduced by {\sc Mader} in~\cite{Mader1988}, generalizing the concept of critical connectivity.
This approach makes it possible to generalize some of the earlier results on foxes from $3$-connected graphs to certain classes of $k$-connected graphs.

Extending the definition above, let us define a {\em $k$-fox} to be a $k$-connected graph admitting a spanning tree without $k$-contractible edges.
For $k \geq 4$, there are graphs $G$ without $k$-contractible edges at all, and every such $G$ is, trivially, a $k$-fox; thus, the question is interesting only under
additional constraints to $G$ which force $k$-contractible edges. Classic constraints are to forbid triangles or to bound the vertex degrees from below:
In~\cite{Thomassen1981} it has been proven that every triangle free $k$-connected graph contains a $k$-contractible edge,
and in~\cite{Egawa1991}, it has been shown that every $k$-connected graph of minimum degree at least $\frac{5k-3}{4}$ must contain a $k$-contractible edge
(unless $G$ is isomorphic to $K_{k+1}$ when $k \leq 3$).
These results do have a common root in terms of generalized criticity~\cite{Mader1988}, and so it is perhaps not surprising that the following new result,
Theorem~\ref{T2}, follows from a statement on special separating sets (Theorem~\ref{T5} in Section~\ref{SGeneral}).

\begin{theorem}
  \label{T2}
  Let $G$ be a $k$-connected graph (except for $K_{k+1}$ if $k \in \{1,2\}$) that is triangle free or of minimum degree at least $\frac{3}{2}k-1$.
  Then every spanning tree of $G$ contains at least two $k$-contractible edges.
\end{theorem}

This implies that $k$-foxes must contain triangles as well as vertices of ``small'' degree.
In order to show that the bound in Theorem~\ref{T2} is best possible, we exhibit $k$-connected graphs of minimum degree $\frac{3}{2}k-\frac{3}{2}$ (and necessarily containing triangles) that admit a spanning tree with no $k$-contractible edge.
For odd $k \geq 3$, take the lexicographic product of any cycle and $K_{(k-1)/2}$ and add an additional vertex plus all edges connecting it to the others.
(So for $k=3$ we get back the wheels.) The resulting graph is $k$-connected and of minimum degree $\frac{3}{2}k-\frac{3}{2}$,
and the spanning star centered at the additional vertex has no $k$-contractible edge. The same construction works, more generally, if instead of a cycle we start
with any {\em critically $2$-connected graph}, that is, a $2$-connected graph $G$ such that for every vertex $x$ the graph $G-x$ is not $2$-connected. ---
However, for DFS trees the situation changes once more:

\begin{theorem}
  \label{T3}
  For $k>3$, every DFS tree of every $k$-connected graph of minimum degree at least $\frac{3}{2}k-\frac{3}{2}$ 
  contains at least two $k$-contractible edges.
\end{theorem}

Observe that the statement of Theorem~\ref{T3} remains true for $k=3$ by Theorem~\ref{T1} unless the graph is one of the three exceptions listed there.

Theorem~\ref{T2} provides a particularly simple proof that every spanning tree of a {\em cubic} $3$-connected graph nonisomorphic to $K_4$ or the prism has
at least two $3$-contractible edges (see Corollary~\ref{C2} in Section~\ref{SGeneral}); however, taking more external knowledge into account we can
improve {\em two} to the following sharp linear bound in terms of $|V(G)|$ (end of Section~\ref{SGeneral}).

\begin{theorem}
  \label{T4}
	Every spanning tree of every cubic $3$-connected graph nonisomorphic to $K_4$ contains at least $\frac{1}{3}|V(G)|-1$ many $3$-contractible edges.
  The bound is sharp, also when restricted to DFS trees.
\end{theorem}

We also show sharpness for Theorem~\ref{T4}. Obtain a graph $G'$ from any cubic $3$-connected graph $G$ by replacing every vertex $x$ with a triangle $\Delta_x$ such that, for every incident edge $e$ of $x$, the end vertex $x$ of $e$ is replaced with a unique vertex of $\Delta_x$. Clearly, $G'$ is cubic and $3$-connected. Let $T$ be a spanning tree of $G$, and let $T'$ be formed by all edges of $T$ together with the edges of a spanning path of each $\Delta_x$. Then $T'$ is a spanning tree of $G'$ with exactly $\frac{1}{3}|V(G')|-1$ many $3$-contractible edges, as no edge in a triangle is $3$-contractible. When restricted to DFS trees, assume in addition that $G$ is Hamiltonian and let $T$ be a Hamiltonian path of $G$. Then the paths of each $\Delta_x$ can be chosen such that $T'$ is a Hamiltonian path of $G'$ and we see that there is no improvement for DFS trees in general.

\section{Contractible edges in spanning trees}
\label{SGeneral}

Let $G$ be a graph and $\mathfrak{T}(G)$ be the set of its smallest separating sets.
For $T \in \mathfrak{T}(G)$, the union of the vertex sets of at least one but not of all components of $G-T$ is called a {\em $T$-fragment}.
Obviously, if $F$ is a $T$-fragment then so is $\overline{F}^G:=V(G) \setminus (T \cup F)$, where the index $G$ is always omitted as it will be clear from the context.
Moreover, $\overline{\overline{F}}=F$. Fragments have the following fundamental property.

\begin{lemma}\emph{\cite{Mader1988}}
  \label{L1}
  Let $B$ be a $T_B$-fragment and $F$ be a $T$-fragment of a graph $G$ such that $B \cap F \not= \emptyset$.
  Then $|B \cap T| \geq |\overline{F} \cap T_B|$, and if equality holds then $B \cap F$ is a $(B \cap T) \cup (F \cap T_B) \cup (T \cap T_B)$-fragment.
\end{lemma}

{\bf Proof.} Let $k:=\kappa(G)$ and observe that $N_G(B \cap F)$ separates $G$ and is a subset of $X:=(B \cap T) \cup (F \cap T_B) \cup (T \cap T_B)$.
Therefore, $k \leq |N_G(B \cap F)| \leq |X|=|B \cap T|+|T_B| - |T_B \cap \overline{F}|=|B \cap T|+k-|T_B \cap \overline{F}|$. Since $k$ cancels on both sides, rearranging the terms yields the desired inequality, and equality implies $N_G(B \cap F)=X$.
\hspace*{\fill}$\Box$

We will not give explicit references to Lemma~\ref{L1}, but mark estimations or conclusions based on it by $\star$;
for example, we write $|F \cap T'| \stargeq |F' \cap T|$ if $F$ is a $T$-fragment and $F'$ is a $T'$-fragment such that $F \cap \overline{F'} \not= \emptyset$
to indicate that the inequality is a straightforward application of Lemma~\ref{L1}.
This convention also applies to the following slightly more complex but standard application of Lemma~\ref{L1}: 
If both $B \cap F$ and $\overline{B} \cap \overline{F}$ are nonempty, then, by Lemma~\ref{L1}, they are both fragments. In many cases, $B$ will be an inclusion minimal fragment with respect to some property, $F$ will be a $T$-fragment such that $T$ contains a vertex from $B$,
and $F \cap B \not= \emptyset$ will have the same property as $B$ (but is no fragment by minimality): In such a scenario, we infer 
$|B \cap T| \geq |\overline{F} \cap T_B|+1$, $|F \cap T_B| \geq |\overline{B} \cap T|+1$, and $\overline{F} \cap \overline{B}=\emptyset$
from Lemma~\ref{L1}, and again refer to it by $\star$, for example, by writing
$|B \cap T| \stargeq |\overline{F} \cap T_B|+1$, $|F \cap T_B| \stargeq |\overline{B} \cap T|+1$, or $\overline{F} \cap \overline{B} \stareq \emptyset$, respectively.

Another fact that will be used throughout is that if $F$ is a $T$-fragment contained in some smallest separating set $T'$, then $T$ contains
a vertex of every $T'$-fragment (as every vertex from $F \subseteq T'$ must have a neighbor in every $T'$-fragment, which can only be in $T$);
since every component of $\overline{F}$ is adjacent to all vertices of $T$, it must contain a vertex from $T'$, so that, in particular, $\overline{F} \cap T' \not= \emptyset$.

Now let us fix a subset $\mathfrak{S}$ of the power set $\mathfrak{P}(V(G))$. We call a $T$-fragment $F$ a {\em $T$-$\mathfrak{S}$-fragment} if
$S \subseteq T$ for some $S \in \mathfrak{S}$. In that case, again, $\overline{F}$ is a $T$-$\mathfrak{S}$-fragment, too;
$F$ is called a {\em $T$-$\mathfrak{S}$-end} if there is no $T'$-$\mathfrak{S}$-fragment properly contained in it,
and $F$ is called a {\em $T$-$\mathfrak{S}$-atom} if there does not exist a $T'$-$\mathfrak{S}$-fragment with fewer than $|F|$ vertices.
Observe that if $F$ is a $T$-fragment then necessarily $T=N_G(F)$, so that $T$ can be reconstructed from $F$; therefore,
one might omit $T$ in the notion, which defines the terms {\em fragment}, {\em $\mathfrak{S}$-fragment}, {\em $\mathfrak{S}$-end}, and {\em $\mathfrak{S}$-atom}.
These definitions and the following theorem are from~\cite{Mader1988} and have their roots back in a 1970 paper by {\sc Watkins}
where it was proven that the degrees of a vertex transitive $k$-connected graph are at most $\frac{3}{2}k-1$~\cite{Watkins1970}.

\begin{theorem}\emph{\cite{Mader1988}}
  \label{TMaderGeneral} 
  Let $G$ be a graph, $\mathfrak{S} \subseteq \mathfrak{P}(V(G))$, and $A$ be a $T_A$-$\mathfrak{S}$-atom of $G$. Suppose that there exists
  an $S \in \mathfrak{S}$ and a $T \in \mathfrak{T}(G)$ such that $S \subseteq T \setminus \overline{A}$ and $T \cap A \not= \emptyset$.
  Then $A \subseteq T$ and $|A| \leq |T \setminus T_A|/2$.
\end{theorem}

A fragment of minimum size is usually called an {\em atom} of $G$. Consequently, for $\mathfrak{S}:=\{\{x\}:\, x \in V(G)\}$, we obtain the following specialization of Theorem~\ref{TMaderGeneral},
which appeared already in~\cite{Mader1971}.

\begin{theorem}\emph{\cite{Mader1971}}
  \label{TMaderSpecial} 
  Let $G$ be a graph and $A$ be a $T_A$-atom of $G$. Suppose that $A \cap T \not= \emptyset$ for some $T \in \mathfrak{T}(G)$.
  Then $A \subseteq T$ and $|A| \leq |T \setminus T_A|/2 \leq \kappa(G)/2$.
\end{theorem}

We start our considerations with the following result.

\begin{theorem}
  \label{T5}
  Let $Q$ be a spanning tree of a 
  graph $G$ of connectivity $k$, set
  $\mathfrak{S}:=\{V(e):\, e \in E(Q)\}$, suppose that all $\mathfrak{S}$-fragments have cardinality at least $\frac{k-1}{2}$,
  and let $B$ be an $\mathfrak{S}$-end. Then $|B|=\frac{k-1}{2}$ (in particular, $k$ is odd) or all edges $e$ from $Q$ with $|V(e) \cap B|=1$ are $k$-contractible.
\end{theorem}

{\bf Proof.}
Let $a:=\frac{k-1}{2}$ and $T_B:=N_G(B)$. Observe that the existence of $B$ implies $k > 1$. Since $Q$ is a spanning tree, there exists an edge $e$ with $|V(e) \cap B|=1$.
If all such edges $e$ are $k$-contractible then we are done.
Otherwise, one such edge $e$ is not $k$-contractible; there exists a $T \in \mathfrak{T}(G)$ with $V(e) \subseteq T$, and we consider a $T$-fragment $F$.
Now $B$ and $F$ are $\mathfrak{S}$-fragments, so that $|B|,|F|,|\overline{B}|,|\overline{F}| \geq a$, and it suffices to prove that $|B| \leq a$.

Observe that $V(e) \cap T_B \not=\emptyset$, so $|T \cap T_B| \geq 1$.
If $B \cap F \not= \emptyset \not= B \cap \overline{F}$, we infer $\overline{B} \starsubseteq T$
and $2|\overline{B}| = 2|\overline{B} \cap T| \starleq |F \cap T_B|-1+|\overline{F} \cap T_B|-1 = |T_B \setminus T|-2 \leq k-3$, and, hence, $|\overline{B}| \leq (k-3)/2 < a$, which is a contradiction.
Suppose that $B \cap F \not= \emptyset$. Then $B \cap \overline{F}=\emptyset$, which implies $\overline{F} \starsubseteq T_B$.
If $\overline{B} \cap F \not = \emptyset$, then $|\overline{B} \cap T| \stargeq |\overline{F} \cap T_B| =|\overline{F}| \geq a$, and otherwise $|\overline{B} \cap T|=|\overline{B}| \geq a$, too.
Hence, $k=|T| = |B \cap T|+|\overline{B} \cap T|+|T_B \cap T| \stargeq (|\overline{F} \cap T_B|+1)+a+1 \geq 2a+2=k+1$, which is a contradiction.
Consequently, $B \cap F=\emptyset$, and, by the same argument, $B \cap \overline{F}=\emptyset$. 
If $\overline{B} \cap F \not= \emptyset$, we infer
$|F \cap T_B| \stargeq |B \cap T|=|B| \geq a$, and otherwise $|F \cap T_B|=|F| \geq a$, too.
Symmetrically, we get $|\overline{F} \cap T_B| \geq a$ and conclude $|T \cap T_B|=1$.
Now $\overline{B} \cap F \not= \emptyset$ implies $|\overline{B} \cap T| \stargeq |\overline{F} \cap T_B| \geq a$,
$\overline{B} \cap \overline{F} \not= \emptyset$ implies $|\overline{B} \cap T| \geq^\star |F \cap T_B| \geq a$, and otherwise $\overline{B} \cap F = \emptyset = \overline{B} \cap \overline{F}$ implies $|\overline{B} \cap T|=|\overline{B}| \geq a$, too. It follows $|B|=|B \cap T| = |T|-|T \cap T_B|-|\overline{B} \cap T| \leq k-1-a=a$, which proves the theorem.
\hspace*{\fill}$\Box$

\begin{corollary}
  \label{C1}
  Let $G \not\cong K_{k+1}$ 
  be a $k$-connected graph in which every fragment has cardinality at least $\frac{k}{2}$.
  Then every spanning tree of $G$ admits at least two $k$-contractible edges.
\end{corollary}

{\bf Proof.} 
Let $Q$ be a spanning tree of $G$. Then $|E(Q)| \geq 2$, since $|V(G)| \geq k+2 \geq 3$. Hence, we may assume that at least one edge in $Q$ is not $k$-contractible. Therefore, $\kappa(G)=k$, and there exists an $\mathfrak{S}$-fragment $C$, where we define $\mathfrak{S}:=\{V(e):\, e \in E(Q)\}$ as before. In particular, $k > 1$. Consider an $\mathfrak{S}$-end $B \subseteq C$.
By Theorem~\ref{T5}, $Q$ contains a $k$-contractible edge $e$ that has precisely one end vertex in $C$.
Likewise, consider an $\mathfrak{S}$-end $B \subseteq \overline{C}$. By applying Theorem~\ref{T5} once more, $Q$ contains a $k$-contractible edge $f$ that has precisely one end vertex in $\overline{C}$.
Clearly, $e \not= f$, which proves the statement.
\hspace*{\fill}$\Box$

For $k \leq 2$, the fragment size condition of Corollary~\ref{C1} is trivially true, and so every spanning tree of every $k$-connected graph nonisomorphic to $K_{k+1}$, $k \leq 2$, admits at least two $2$-contractible edges. For general $k$, Corollary~\ref{C1} implies Theorem~\ref{T2} as follows, and the examples beneath the latter in the introduction show also that the bound on the fragment size in Corollary~\ref{C1} cannot be improved.

{\bf Proof of Theorem~\ref{T2}.}
As argued above, the statement holds for $k \leq 2$, so let $k \geq 3$. Then $G \neq K_{k+1}$, since $G$ is triangle free or of minimum degree at least $\frac{3}{2}k-1$. If the minimum degree condition is satisfied, every fragment has cardinality at least $\frac{k}{2}$, and applying Corollary~\ref{C1} gives the claim. If $G$ is triangle free, let $Q$ be a spanning tree of $G$ and let $\mathfrak{S}:=\{V(e):\, e \in E(Q)\}$. As in the proof of Corollary~\ref{C1}, we may assume that at least one edge in $Q$ is not $k$-contractible, which implies $\kappa(G)=k$ and the existence of a $\mathfrak{S}$-fragment $C$. Since $G$ is triangle free, every $\mathfrak{S}$-fragment contains two adjacent vertices, and considering the neighborhood of the two vertices of degree at least $k$ each implies that every $\mathfrak{S}$-fragment has in fact cardinality at least $k$.
By applying Theorem~\ref{T5} twice, as in the previous proof, we find two $k$-contractible edges $e \not= f$ in $Q$ with end vertices in $C$ and $\overline{C}$, respectively.
\hspace*{\fill}$\Box$

As promised in the introduction we derive the following result from Theorem~\ref{T2}.

\begin{corollary}
  \label{C2}
  Every spanning tree of every cubic $3$-connected graph nonisomorphic to $K_4$ or the prism $K_3 \times K_2$ contains at least two $3$-contractible edges.
\end{corollary}

{\bf Proof.} We use induction on the number of vertices.
The induction starts for $K_4$, so suppose that $G$ is a cubic $3$-connected graph on at least six vertices, and $Q$ is a spanning tree of $G$.
We may assume that $G$ contains a triangle $\Delta$, for otherwise the statement follows from Theorem~\ref{T2}, and we may assume that $G$ is not the prism.
The edge neighborhood of any such triangle forms a matching of three $3$-contractible edges, at least one of which belongs to $Q$;
in fact, we may assume that exactly one edge of the edge neighborhood belongs to $Q$, for otherwise the statement is proven.
So suppose that $e$ is the only edge from $Q$ in the edge neighborhood of $\Delta$.
If there was another triangle $\Delta'$ then, consequently, $e$ is the only edge from $Q$ in the edge neighborhood of $\Delta'$, too, implying that $V(G)=V(\Delta) \cup V(\Delta')$,
so that $G$ is the prism $K_3 \times K_2$, which is a contradiction. So we may assume that $\Delta$ is the only triangle in $G$.
The graph $G/\Delta$ obtained from $G$ by identifying the three vertices of $\Delta$ and simplifying
is not $K_4$, as $G$ is not the prism, and $G/\Delta$ is not the prism, as $\Delta$ is the only triangle in $G$ and the prism has two vertex-disjoint triangles.
Clearly, $G/\Delta$ is cubic and $Q/\Delta$ is a spanning tree of $G$. Since no smallest separating set of $G$ contains two vertices of $\Delta$, $G/\Delta$ is also $3$-connected.
Hence, by induction, $Q/\Delta$ contains two $3$-contractible edges of $G/\Delta$, and the two edges corresponding to these in $G$ are $3$-contractible in $G$ as one checks readily.\qed

By using a powerful result on $3$-contractible edges in $3$-connected graphs from the literature we can improve Corollary~\ref{C2} to Theorem~\ref{T4}
(where the lower bound to the number of $3$-contractible edges is sharp).

{\bf Proof of Theorem~\ref{T4}.}
From Lemma~3 and Lemma~4 in~\cite{Kriesell2000}, we get Theorem~3 in~\cite{Kriesell2002},
which implies, together with Theorem~12 from~\cite{Kriesell2002},
that {\em every vertex in a $3$-connected graph is either contained in a triangle or on at least two $3$-contractible edges}.
We first show that this implies for any $3$-connected {\em cubic} graph $G$ that its subgraph $H$ on $V(G)$ formed by the non-$3$-contractible edges of $G$ is a clique factor (that is, all components of $H$ are isolated vertices, or single edges, or triangles --- or $K_4$ in case that $G$ is $K_4$): Suppose that $G$ is not $K_4$. If $x$ has degree at least $2$ in $H$, then $x$ is on a triangle $\Delta$ in $G$ by the result mentioned above and this triangle is also in $H$, whereas the edges from its edge neighborhood are not (so $x$ is in a triangle component of $H$). Now every spanning tree $Q$ contains at most $2$ edges
from every triangle in $H$, so that it contains at most $\frac{2}{3}|V(G)|$ edges from $H$. Therefore, $Q$ contains at least $|V(G)|-1-\frac{2}{3}|V(G)|=\frac{1}{3}|V(G)|-1$ $3$-contractible edges (unless $G$ is $K_4$).\hspace*{\fill}$\Box$

\section{Contractible edges in DFS trees}

Again, we observe that the spanning tree in the sharpness example of Corollary~\ref{C1} and Theorem~\ref{T2} is far from being a DFS tree.
The following theorem will provide more insight into the distribution of $k$-contractible edges in spanning trees of graphs where
the fragment lower bound $\frac{k-1}{2}$ from Theorem~\ref{T5} is sharp (as opposed to the bound $\frac{k}{2}$ that is sharp in Corollary~\ref{C1}),
and leads to a proof of Theorem~\ref{T3} and, in the next section, of Theorem~\ref{T1}.

\begin{theorem}
  \label{T6}
  Let $Q$ be a spanning tree of a 
  graph $G$ of connectivity $k$, set
  $\mathfrak{S}:=\{V(e):\, e \in E(Q)\}$, suppose that all $\mathfrak{S}$-fragments have cardinality at least $\frac{k-1}{2}$, set
  $\mathfrak{R}:=\{V(e):\, e \in E(Q)$, $|V(e) \cap A|=1$ for some $\mathfrak{S}$-end $A\}$, and let $B$ be an $\mathfrak{R}$-end.
  Then $|B|=\frac{k-1}{2}$, or $Q$ contains a $k$-contractible edge $e$ with 
  at least
  one endvertex in $B$,
  or $N_G(B)$ contains an $\mathfrak{R}$-fragment of cardinality $\frac{k-1}{2}$ such that all edges from $Q$ having exactly one vertex in common with it are $k$-contractible.
\end{theorem}

{\bf Proof.}
Let $a:=\frac{k-1}{2}$. 
Let us call an edge from $Q$ {\em green} if it is not $k$-contractible.

{\bf Claim~1.} Suppose that $e$ is a green edge and $A$ is an $\mathfrak{S}$-end with $|V(e) \cap A|=1$.
Then $|A|=a$ (so that $A$ is an $\mathfrak{S}$-atom, and $a$ is an integer), 
and $A \subseteq T$ for every $T \in \mathfrak{T}(G)$ with $V(e) \subseteq T$ (and there exists such a $T$).

We get $|A|=a$ immediately from Theorem~\ref{T5}. Since $e$ is green, there exists a $T \in \mathfrak{T}(G)$ such that $V(e) \subseteq T$, and
for every such $T$ we know $V(e) \subseteq T \setminus \overline{A}$ and $A \cap T \not= \emptyset$, so that
$A \subseteq T$ follows from Theorem~\ref{TMaderGeneral}, proving Claim~1.

Now let us call a green edge $e$ {\em red} if $|V(e) \cap A|=1$ for some $\mathfrak{S}$-end $A$ (see Figure~\ref{fig:LargeDFSTreeWith2ContractibleEdges} for an example of this coloring in the special case $k=3$). 
Let $T_B:=N_G(B)$. By definition of $\mathfrak{R}$ and Claim~1, there
exists a red edge $e'$ and an ${\mathfrak S}$-end $A'$ with $|V(e') \cap A'|=1$ and $A' \cup V(e') \subseteq T_B$.
Since $B$ is an ${\mathfrak S}$-fragment (every $\mathfrak{R}$-fragment is an ${\mathfrak S}$-fragment by definition), it must contain an $\mathfrak{S}$-end $A$. There exists an edge $e$ from $Q$
with $|V(e) \cap A|=1$. If $e$ is $k$-contractible, we are done; thus we may assume that $e$ is green, so that $e$ is even red.
By definition and Claim~1, there exists a $T \in {\mathfrak T}(G)$ such that $A \cup V(e) \subseteq T$. We now consider a $T$-fragment $F$,
which is, in fact, an $\mathfrak{R}$-fragment, and analyze the possible ways $F,\overline{F}$ meet $B,\overline{B}$ (for example, $B=A=\{v\}$, $e'=ya$ and $F=A'=\{a\}$ in Figure~\ref{fig:LargeDFSTreeWith2ContractibleEdges}). 

Again we may rule out that $B \cap F \not= \emptyset \not= B \cap \overline{F}$, as this would imply $\overline{B} \subseteq^\star T$ and $2|\overline{B}|=2|\overline{B} \cap T| \starleq |\overline{F} \cap T_B|-1 + |F \cap T_B|-1 \leq k-2$, and hence $|\overline{B}| \leq \frac{k-2}{2}<a$, which gives a contradiction.

\begin{figure}[htb]
	\centering
	\subfloat{
		\makebox[4.5cm]{
		\includegraphics[scale=0.7]{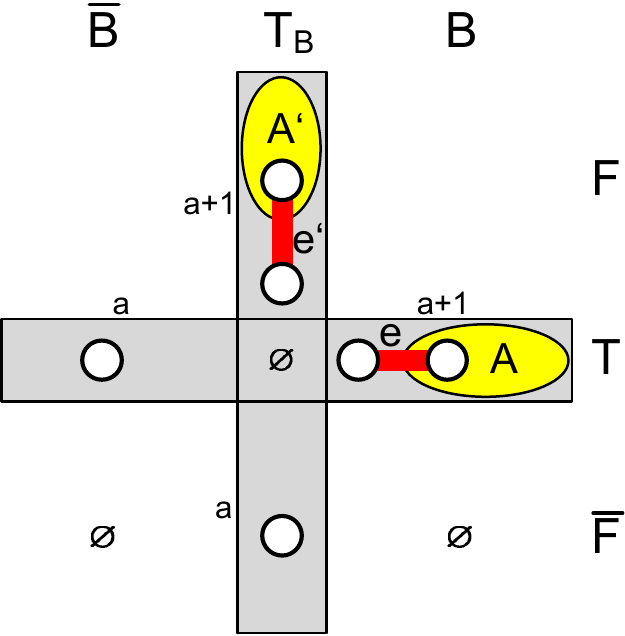}
		\label{fig:SEnd3}
		}
	}
	\caption{The structure of the $\mathfrak{R}$-fragments $B$ and $F$.}
	\label{fig:Structure}
\end{figure}

Now assume that exactly one of $B \cap F$ and $B \cap \overline{F}$ is nonempty, say, by symmetry of $F$ and $\overline{F}$, $B \cap F \not= \emptyset$ (see Figure~\ref{fig:Structure}). It follows $B \cap \overline{F}=\emptyset \stareq \overline{B} \cap \overline{F}$.
Consequently, $|\overline{F} \cap T_B|=|\overline{F}| \geq a$ and, hence, $|B \cap T| \stargeq |\overline{F} \cap T_B|+1 \geq a+1$.
If $\overline{B} \cap F \not= \emptyset$, then $|\overline{B} \cap T| \stargeq |\overline{F} \cap T_B| \geq a$, and if otherwise $\overline{B} \cap F= \emptyset$, then $|\overline{B} \cap T|=|\overline{B}| \geq a$, too. Now $k=|T|=|B \cap T|+|\overline{B} \cap T|+|T_B \cap T| \geq (a+1) + a + 0 = k$, so we get equality summand-wise,
that is, $|B \cap T|=a+1$, $|\overline{B} \cap T|=a$, and $T \cap T_B=\emptyset$. This implies $a \leq |\overline{F}| \leq^\star |B \cap T|-1=a$, so
that $|\overline{F}|=a$ (in particular, $\overline{F}$ is an $\mathfrak{R}$-atom) and $B \cap T=A \cup V(e)$. 
Moreover, since $A' \cup V(e') \subseteq T_B$ and $|A' \cup V(e')| \geq a+1 > |\overline{F}|$, we conclude $F \cap T_B = A' \cup V(e')$, as $G[A' \cup V(e')]$ is connected.

In the remaining case $B \cap F = \emptyset = B \cap \overline{F}$, i.e.\ $B \subseteq T$, we essentially obtain the same picture (Figure~\ref{fig:Structure}):
If $F \cap \overline{B} \not = \emptyset$, then $|F \cap T_B| \stargeq |B \cap T|=|B| \geq a$, and if otherwise $F \cap \overline{B}=\emptyset$, then $|F \cap T_B|=|F| \geq a$, too; symmetrically, $|\overline{F} \cap T_B| \geq a$. Then $F \cap \overline{B} \not= \emptyset$ implies $|\overline{B} \cap T| \stargeq |\overline{F} \cap T_B| \geq a$,
$\overline{F} \cap \overline{B} \not= \emptyset$ implies $|\overline{B} \cap T| \geq |F \cap T_B| \geq a$,
and the remaining case $F \cap \overline{B}=\emptyset=\overline{F} \cap \overline{B}$ implies $|\overline{B} \cap T|=|\overline{B}| \geq a$, too. Now if $|B|=a$, the theorem is proved, so let $|B|>a$. Then $k=|T|=|B|+|T \cap T_B|+|\overline{B} \cap T| \geq (a+1)+0+a=k$, so that equality holds summand-wise, implying
$T \cap T_B=\emptyset$, $B=B \cap T=A \cup V(e)$, $|B|=a+1$, and $|\overline{B} \cap T|=a$.
By symmetry of $F$ and $\overline{F}$, we may assume that $A' \cup V(e') \subseteq F$ and hence $|F \cap T_B|=a+1$ and $|\overline{F} \cap T_B|=a$. This implies $\overline{B} \cap \overline{F}=\emptyset$, as otherwise $k=|F \cap T_B|+|\overline{F} \cap T_B| \stargeq |A' \cup V(e')|+|B \cap T|=(a+1)+(a+1)>k$ gives a contradiction.

Regardless whether $B \cap F$ is empty or not, we conclude $T \cap T_B=\emptyset$, $B \cap T=A \cup V(e)$ with cardinality $a+1$, $F \cap T_B=A' \cup V(e')$ with cardinality $a+1$, and $\overline{F} = \overline{F} \cap T_B$ is an $\mathfrak{R}$-atom with cardinality $a$ and thus also an $\mathfrak{S}$-atom. We proceed with the general argument.

If all edges from $Q$ that connect $\overline{F}$ to $T$ are $k$-contractible, we are done.
So there exists a green edge $xy$ with $x \in \overline{F}$ and $y \in T$ and a $T' \in \mathfrak{T}(G)$ with $\{x,y\} \subseteq T'$,
and by Claim~1 (applied to $xy$ and $\overline{F}$) we obtain $\overline{F} \subseteq T'$.
We discuss the possible locations of $y$ and will show that all are impossible, which proves the theorem.
Observe that $A$ must have at least $a=|A|$ neighbors in $\overline{F}$,
for otherwise $\overline{F} \setminus N_G(A)$ would be a nonempty set with less than $k$ neighbors in $G$, which contradicts the fact that $G$ is $k$-connected.
It follows $\overline{F} \subseteq N_G(A)$.

If $y$ would be the vertex in $V(e) \setminus A$, then $xy$ would be a red edge with its endvertices in $N_G(A)$,
certifying that $A$ is an $\mathfrak{R}$-fragment properly contained in the $\mathfrak{R}$-end $B$,
which is absurd. 
If $y$ would be some vertex in $A$, then $xy$ would be a red edge, and by Claim~1 (applied to $xy$), we get $A \subseteq T'$.
As $B$ is an $\mathfrak{R}$-end, $T'$ and $T_B$ \emph{cross} (i.e.\ $T' \cap B \neq \emptyset \neq T' \cap \overline{B}$
), and $T'$ contains at least one vertex from $\overline{B}$; on the other hand, $T'$ contains
$a+a=k-1$ vertices from $A \cup \overline{F} \subseteq B \cup T_B$, so that $T'$ contains exactly one vertex from $\overline{B}$.
Consider a $T'$-fragment $F'$ that contains $T_B \cap F$ (this exists, as $G[T_B \cap F]$ is connected and does not intersect $T'$).
Then $T_B$ must intersect $\overline{F'}$, which contradicts $T_B \setminus F = \overline{F} \subseteq T'$.

It follows that, necessarily, $y \in \overline{B}$. Suppose to the contrary that $T'$ separates $A \cup V(e)$, that is, there exists a $T'$-fragment $F'$ such
that $F' \cap (A \cup V(e)) \not= \emptyset$ and $\overline{F'} \cap (A \cup V(e)) \not= \emptyset$. Then $T'$ must contain a vertex $z$ from $A$, and $T_B$ and $T'$ cross,
so that $T'$ separates $A' \cup V(e') \supseteq T_B \setminus T'$, too. Thus, $T'$ contains a vertex $z'$ from $A'$. Since $T'$ cannot separate the end vertices of $V(e')$, we know that
$V(e') \cap F' = \emptyset$ or $V(e') \cap \overline{F'}=\emptyset$. Without loss of generality we may suppose that $V(e') \cap \overline{F'}=\emptyset$.
If $F' \cap \overline{B} \not= \emptyset$, then $F' \cap \overline{B}$ would be an $S :=(F' \cap T_B) \cup (T' \cap T_B) \cup (T' \cap \overline{B})$-fragment, as $\overline{F'} \cap B \not= \emptyset$.
Since $e'$ is a green edge with exactly one end vertex from the $\mathfrak{S}$-atom $A'$ and contained in $S$ we know from Claim~1 that $A' \subseteq S$, contradiction. 
If $\overline{F'} \cap \overline{B} \not= \emptyset$, then $F' \cap B$ would be an $S := (F' \cap T_B) \cup (T' \cap T_B) \cup (T' \cap B)$-fragment containing $V(e')$, as $F' \cap B \not= \emptyset$. This implies again $A' \subseteq S$, contradiction. 
Therefore, $\overline{B} \subseteq T'$. It follows $k=|T'|=|\overline{B} \cap T'|+|T_B \cap T'|+|B \cap T'| \geq |\overline{B}|+|\overline{F} \cup \{z'\}|+|\{z\}| \geq a+(a+1)+1>k$, contradiction.

Hence we have to assume that $T'$ does not separate $A \cup V(e)$. Consequently, there exists a $T'$-fragment $F'$ such that $F' \cap (A \cup V(e))=\emptyset$.
If $F' \cap F \not= \emptyset$, then $|F' \cap T| \stargeq |\overline{F} \cap T'|=|\overline{F}| \geq a$, and if otherwise $F' \cap F = \emptyset$, then $|F' \cap T|=|F'| \geq a$, too. On the other hand, $|F' \cap T| \leq k - |A \cup V(e)|-|\{y\}| = k-(a+1)-1=a-1$, contradiction.
\hspace*{\fill}$\Box$

We now prove a condition that guarantees two $k$-contractible edges in any spanning tree.

\begin{lemma}
  \label{L2}
  Let $Q$ be a spanning tree of a noncomplete graph $G$ of connectivity $k$, set
  $\mathfrak{S}:=\{V(e):\, e \in E(Q)\}$, suppose that all $\mathfrak{S}$-fragments have cardinality at least $\frac{k-1}{2}$, set
  $\mathfrak{R}:=\{V(e):\, e \in E(Q)$, $|V(e) \cap A|=1$ for some $\mathfrak{S}$-end $A\}$.
  Then $Q$ contains two $k$-contractible edges or there is an $\mathfrak{R}$-end of cardinality $\frac{k-1}{2}$ such that $Q$ contains no $k$-contractible edge with an endvertex from it.
\end{lemma}

{\bf Proof.} 
Call a fragment $C$ {\em small} if $|C|=\frac{k-1}{2}$ and {\em big} otherwise; call $C$ {\em good} if $Q$ contains a $k$-contractible edge that has at least one endvertex in $C$ and {\em bad} otherwise. Call $C$ {\em very good} if all edges from $Q$ having exactly one endvertex in $C$ are $k$-contractible.
In this language, Theorem~\ref{T6} tells us that an $\mathfrak{R}$-end $B$ is small, or good, or $N_G(B)$ contains a small, very good $\mathfrak{R}$-fragment.

We may assume that $Q$ contains at least one non-$k$-contractible edge, so that there exists an $\mathfrak{S}$-fragment, and, hence, an $\mathfrak{S}$-end, say, $A$.
Take an $\mathfrak{S}$-end $A' \subseteq \overline{A}$. There exist edges $e,e' \in E(Q)$ with $|V(e) \cap A|=1$ and $|V(e') \cap A'|=1$; clearly, $e \not= e'$,
and we are done if both $e,e'$ are $k$-contractible. In what remains we thus may assume that there {\em exists} an $\mathfrak{R}$-fragment.

Suppose that there exists a very good $\mathfrak{R}$-fragment $C$.
Consider an $\mathfrak{R}$-end $B \subseteq \overline{C}$. If $B$ is good, then there exists a $k$-contractible edge $e$ having a vertex in common with $B$
and another one having a vertex in common with $C$, which proves the lemma. Hence $B$ is bad. If $B$ is small, then the lemma is proved again,
so we may assume that $B$ is big. By Theorem~\ref{T6}, $N_G(B)$ contains a (small) very good $\mathfrak{R}$-fragment $D$. As $C$ and $D$ are disjoint and their
union is not equal to $V(G)$, $Q$ contains two distinct $k$-contractible edges that are incident with vertices from $C$ or $D$, which proves the lemma.

Therefore, we may assume that there are no very good $\mathfrak{R}$-fragments. Consequently, by Theorem~\ref{T6}, every big $\mathfrak{R}$-end is good.
Moreover, we may assume that every small $\mathfrak{R}$-end is good, for otherwise we are done. Hence, {\em every} $\mathfrak{R}$-end $B$ is good, and, for any $\mathfrak{R}$-end $C$ contained in $\overline{B}$, $Q$ contains two distinct $k$-contractible edges that have an endvertex in $B$ and $C$, respectively, which gives the lemma.
\hspace*{\fill}$\Box$.

We now specialize to DFS trees. A {\em DFS tree} of some graph $G$ is a spanning tree $Q$ with a prescribed {\em root vertex} $r$ such that
for every vertex $x$, any two $x$-branches are nonadjacent in $G$, where an {\em $x$-branch} is the vertex set of any component of $Q-x$ not containing $r$.
Now we are prepared to prove the following theorem, which implies Theorem~\ref{T3} immediately, \
as, in a graph of connectivity $k$ and of minimum degree at least $\frac{3}{2}k-\frac{3}{2}$, every fragment has cardinality at least $\frac{k-1}{2}$.

\begin{theorem}
  \label{T7}
  Let $Q$ be a DFS tree of a noncomplete graph $G$ of connectivity $k>3$, and set $\mathfrak{S}:=\{V(e):\, e \in E(Q)\}$.
  (i) If all $\mathfrak{S}$-fragments have cardinality at least $\frac{k-1}{2}$, then $Q$ contains at least one $k$-contractible edge.
  (ii) If all fragments have cardinality at least $\frac{k-1}{2}$, then $Q$ contains at least two $k$-contractible edges.
\end{theorem}

{\bf Proof.}
First, let us assume that all $\mathfrak{S}$-fragments have cardinality at least $\frac{k-1}{2}$.
Set $\mathfrak{R}:=\{V(e):\, e \in E(Q)$, $|V(e) \cap A|=1$ for some $\mathfrak{S}$-end $A\}$ and observe that $\frac{k-1}{2}>1$. 
By Lemma~\ref{L2}, we may assume that there exists a $T_B$-$\mathfrak{R}$-end $B$ with $|B|=\frac{k-1}{2}$ such that $Q$ contains no $k$-contractible edge with an endvertex from $B$.
There exists an edge $e \in E(Q)$ and a $T_A$-$\mathfrak{S}$-end $A$ such that $|V(e) \cap A|=1$ and $V(e) \subseteq T_B$.
By Claim~1 in the proof of Theorem~\ref{T6}, we see that $|A|=\frac{k-1}{2}$ and $A \subseteq T_B$. Let $x$ be the vertex in $V(e) \setminus A$.

Observe that $B$ is an $\mathfrak{S}$-end (as it is even an $\mathfrak{S}$-atom) and consider any edge $f \in E(Q)$ with $|V(f) \cap B|=1$ (there exists at least one such edge).
Let $y$ be the vertex in $V(f) \setminus B$. As $f$ is not $k$-contractible, there exists a $T \in \mathfrak{T}(G)$ with $V(f) \subseteq T$.
According to Claim~1 in the proof of Theorem~\ref{T6}, $B \subseteq T$. All vertices from $B$ are neighbors of $A$
(for otherwise $(T_B-A) \cup (N_G(A) \cap B)$ would be a separating vertex set of $G$ with less than $k$ vertices), that is, $B \subseteq T_A$.

\begin{figure}[h!tb]
	\centering
	\subfloat[Claim~(i)]{
		\makebox[4.5cm]{
		\includegraphics[scale=0.8]{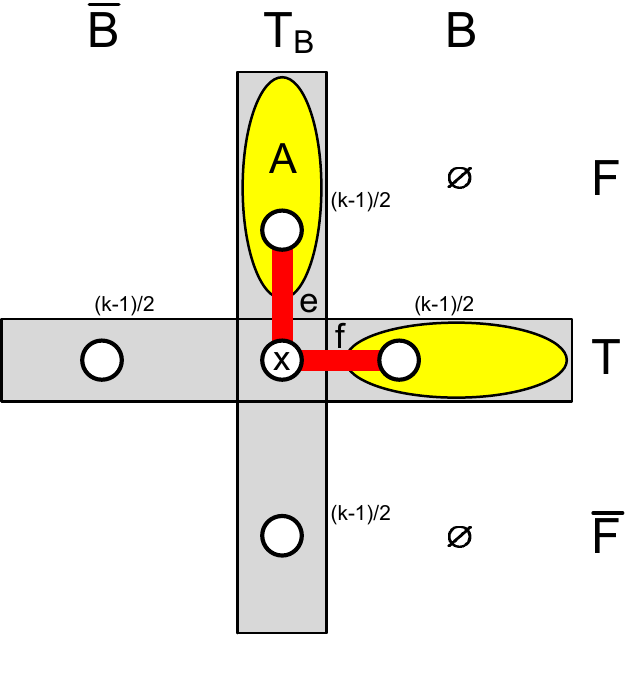}
		\label{fig:SEnd4}
		}
	}
	\hspace{2.5cm}
	\subfloat[Claim~(ii)]{
		\makebox[4.5cm]{
		\includegraphics[scale=0.8]{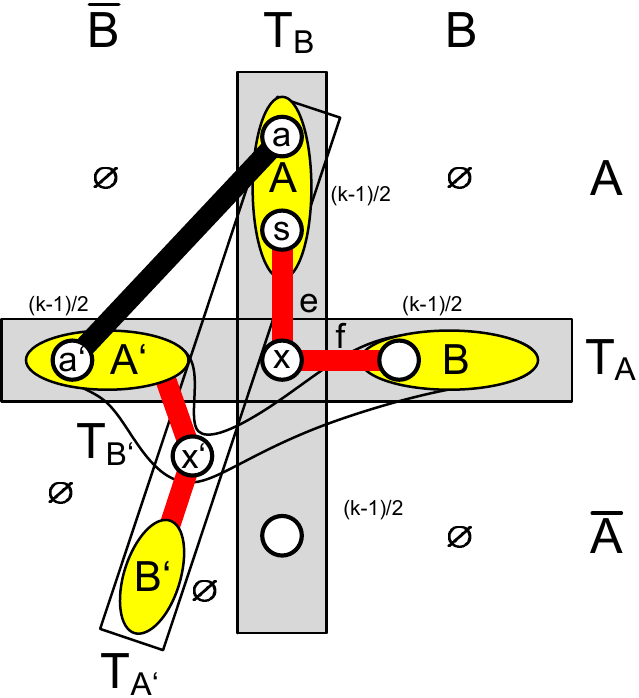}
		\label{fig:SEnd5}
		}
	}
	\caption{The structure of the $\mathfrak{R}$-fragments $B$ and $F$ in Theorem~\ref{T7}.}
	\label{fig:DFSfragments}
\end{figure}

Let $F$ be a $T$-fragment. Then $F \cap \overline{B} \not= \emptyset$ implies $|F \cap T_B| \stargeq |B \cap T|=\frac{k-1}{2}$ and $F \cap \overline{B}=\emptyset$
implies $|F \cap T_B|=|F| \geq \frac{k-1}{2}$, too. As the same holds for $\overline{F}$ instead of $F$, we deduce $|F \cap T_B|=|\overline{F} \cap T_B|=\frac{k-1}{2}$
and $T \cap T_B=\{y\}$. It follows that $A \subseteq F$ or $A \subseteq \overline{F}$, for otherwise $T$ contains a vertex from $A$, as $G[A]$ is connected;
this vertex can only be $y$, so that applying Claim~1 of the proof of Theorem~\ref{T6} on $f$ and $A$ (or alternatively using Theorem~\ref{TMaderGeneral} with $S=V(f)$) implies $A \subseteq T$, which gives the contradiction $|T \cap T_B| \geq |A|>1$. It follows that $A=F \cap T_B$ or $A=\overline{F} \cap T_B$ and $x \in T \cap T_B$, so that $x=y$
(see Figure~\ref{fig:SEnd4}).

Since $f \in E(Q)$ has been chosen arbitrarily from the edge neighborhood of $B$, we see that $x$ is the only neighbor of $B$ in $Q$.
In particular, $A$ is an $\mathfrak{R}$-end (certified by any $f$ as above), since $B \cup \{x\} \subseteq T_A$. Let us take over the notion of small, big, good, bad, and very good fragments from the previous lemma; then $A$ and $B$ are small and $B$ is bad. If $A$ was bad, too, then by symmetry of $A$ and $B$ we see that $x=y$ is the only neighbor of $A$ in $Q$. Since $G[A]$ and $G[B]$ are connected, each of them contains at most one $x$-branch of $Q$; if each of them contains an $x$-branch, these $x$-branches span $G[A]$ and $G[B]$ and, as every vertex in $A$ has a neighbor in $B$ and vice versa, the two $x$-branches are adjacent, which is a contradiction. Hence, one of $G[A],G[B]$ does not contain an $x$-branch, and, thus, contains the root of $Q$, whereas the other one is spanned by a single $x$-branch. By symmetry, let us assume that $B$ contains the root, and consider any vertex $z \in \overline{B} \cap T_A$. Again, $z$ and the vertices in $A$ are in different $x$-branches, but $z$ must have at least one neighbor in $A$, contradiction.

Therefore, we may assume from now on that $A$ is good. Hence, $Q$ contains a $k$-contractible edge $aa'$ with $a \in A$ and $a' \in \overline{B} \cap T_A$ (see also Figures~\ref{fig:6vertices1ContractibleEdgeBW2} and~\ref{fig:6vertices1ContractibleEdgeBW}). This proves~(i) of the statement. 

Now let us assume, in addition, that {\em all} fragments have cardinality at least $\frac{k-1}{2}$,
and (reductio ad absurdum) that $aa'$ is the unique $k$-contractible edge of $Q$.
Observe that all small fragments of $G$ are atoms, so that, by Theorem~\ref{TMaderSpecial}, they are either subsets of or disjoint from any smallest separating set.

Let us apply Theorem~\ref{T6} to an arbitrary $\mathfrak{R}$-end $B' \subseteq \overline{A}$ (see Figure~\ref{fig:SEnd5}).
If $B'$ was good, then $Q$ would contain a $k$-contractible edge that has an endvertex in $B'$, and is thus distinct from $aa'$, contradiction.
Hence $B'$ is bad. 
If $B'$ was big, then $N_G(B')$ would contain a small, very good $\mathfrak{R}$-fragment $C$ by Theorem~\ref{T6}.
Since $aa'$ is the unique $k$-contractible edge in $Q$, we see that $a' \in C$, and, as $C$ is an atom and $a' \in T_A$, we infer $C \subseteq T_A$ by Theorem~\ref{TMaderSpecial}.
Since $V(e)$ is not $k$-contractible, $x \not \in C$, which implies that $C$ contains no vertex from $B$ either. As $B \cup \{x\} \subseteq T_A$, it follows $C=T_A \cap \overline{B}$.
Consequently, $x$ is the only neighbor of $\overline{A}$ in $Q$, as $\overline{A}$ is bad, $C$ is very good and $x$ is the only neighbor of $B$ in $Q$.
Now it is not possible to locate the root vertex of our DFS tree: 
It cannot be in $\overline{A} \cup \{x\}$, because then we find adjacent vertices from $A$ and $B$ in distinct $x$-branches,
it cannot be in $A \cup C$, because then we find adjacent vertices from $B$ and $\overline{A} \cap T_B$ in distinct $x$-branches, and
it cannot be in $B$, because then we find adjacent vertices from $C$ and $N_G(C) \cap \overline{A}$ in distinct $x$-branches --- contradiction.

Therefore, $B'$ is a small, bad $\mathfrak{R}$-end, just as $B$.
We may infer --- just as before for $B$ --- that $T_{B'}:=N_G(B')$ contains a small, good $T_{A'}$-$\mathfrak{R}$-end $A'$.
Since $A' \subseteq T_{B'} \subseteq T_A \cup \overline{A}$ we see that $A$ and $A'$ are disjoint. As $A'$ is good and $aa'$ is the only $k$-contractible edge,
it follows $a' \in A'$, and, as $A'$ is an atom having a vertex in common with $T_A$, we see that $A' \subseteq T_A$ by Theorem~\ref{TMaderSpecial}.
Since $A'$ contains $a' \notin T_B$, we infer for the same reason that $A'$ does not contain $x$, and, as $G[A']$ is connected and $a' \in \overline{B}$, $A'$ does not contain any vertex from $B$. 
It follows $A'=\overline{B} \cap T_A$.
Moreover, we get, as for $A,B$ above, a unique vertex $x'$ in $T_{A'} \cap T_{B'}$ and infer that all non-$k$-contractible edges from $Q$ with exactly one
vertex from $B'$ are incident with $x'$.
 
We claim that $aa'$ is the only edge from $Q$ that connects $A$ and $A'$.
Assume to the contrary there was another one, say, $zz' \in E(Q) \setminus \{aa'\}$, with $z \in A$ and $z' \in A'$.
Since $zz'$ is not $k$-contractible, there exists a $T' \in \mathfrak{T}(G)$ with $z,z' \in T'$, and from Theorem~\ref{TMaderSpecial} we get $A \cup A' \subseteq T'$.
Since $A \subseteq T'$, $T'$ contains at least one vertex from $\overline{A}$ 
and thus consists of the $k-1$ vertices from $A \cup A'$ and another vertex from $\overline{A}$;
but then $T_A \setminus T' = B \cup \{x\}$ induces a connected subgraph of $G$, so $T'$ does not separate $T_A$, contradiction.

Observe that (when the position of the root in $Q$ is neglected), the situation is symmetric in $A,B,a,x$ and $A',B',a',x'$.
We have seen that $N_Q(B)=\{x\}$, $N_Q(A)=\{x,a'\}$, $N_Q(A')=\{x',a\}$, and $N_Q(B')=\{x'\}$. It follows that $N_Q(B \cup A \cup A' \cup B')=\{x,x'\}$.
Let us again analyze the position of the root vertex $r$ of our DFS tree. 
Then the $r,x$-path in $Q$ must enter $x$ using an edge incident to $A$ or $B$~(*), for otherwise $e \in E(Q)$ implies that
the vertex $s$ from $V(e) \setminus \{x\}$ and any vertex in $N_G(s) \cap B$ are in distinct adjacent $x$-branches (contradiction).
Likewise, the $r,x'$-path in $Q$ must enter $x'$ using an edge connecting to $A'$ or $B'$~(*').

Setting $X:=B \cup A \cup A' \cup B'$ we thus see that $r \in X$, for otherwise the second last vertex of every $r,X$-path in $Q$ is either $x$ or $x'$
and its $r,x$- or $r,x'$-subpath violates (*) or (*'), respectively. Furthermore, $x \not= x'$, for otherwise $r \in B$ violates~(*'), $r \in A$ violates~(*') if the $r,x'$-path in $Q$ does not use the edge $aa'$ and $r \in A$ violates~(*) if it does, and we get symmetric violations for $r \in B'$ and $r \in A'$, respectively.
Now we can deduce that $R:=Q[X \cup \{x,x'\}]$ is connected, that is, a subtree of $Q$. Suppose, to the contrary, that $R$ contains a $z$ such that there
is no $r,z$-path in $R$. The $r,z$-path in $Q$ therefore has to use vertices from $V(G) \setminus V(R)$, and, therefore, both $x$ and $x'$.
If $x$ is used last (that is, $x'$ is on the $r,x$-path), then (*) is violated, and otherwise (*) is.

By symmetry of $A,B$ and $A',B'$, we may assume that $r$ is not in $B$. Since $|T_B \cap \overline{A}|=\frac{k-1}{2} > 1$, there
exists a vertex $t \in (T_B \cap \overline{A}) \setminus \{x'\}$; $t$ has a neighbor $z \in B$, and $t,z$ cannot be in different $v$-branches for any vertex $v$,
so that either $z$ is on the $r,t$-path in $Q$ or $t$ is on the $r,z$-path in $Q$. The first option cannot occur, as $B$ has only one neighbor in $Q$,
and the second one implies that $t$, like {\em all} vertices from the $r,z$-path, is in $V(R)$. Since $t \not\in \{x,x'\}$ and $t \not\in B \cup A \cup A'$,
we deduce $t \in B'$, in other words, $T_B \cap B' \not=\emptyset$. 
By Theorem~\ref{TMaderSpecial}, $B' \subseteq T_B$, which determines $T_B=A \cup \{x\} \cup B'$, $T_A=A' \cup \{x\} \cup B$,
$T_{A'}=A \cup \{x'\} \cup B'$,
and $T_{B'}=A' \cup \{x'\} \cup B$.
The set $Y:=\overline{A} \cap \overline{B} \cap \overline{A'} \cap \overline{B'}$ has all its neighbors in
$T_A \cup T_B \cup T_{A'} \cup T_{B'} \setminus (A \cup B \cup A' \cup B')=\{x,x'\}$, which implies
that $Y$ is empty. Therefore, $\overline{A} \cap \overline{B}=\{x'\}$,
so that $N_G(x')=A' \cup B' \cup \{x\}$.
Hence $\{x'\}$ is a fragment of cardinality $1$, contradiction.
\hspace*{\fill}$\Box$

\section{Contractible edges in DFS trees of 3-con\-nec\-ted graphs}
\label{Sthree}

For proving Theorem~\ref{T1}, observe that the fragment size conditions from Theorem~\ref{T7} are vacuously true in the case that $k=3$. However, the conclusion in~(ii) does not hold for that case, as shown for its direct implication Theorem~\ref{T3} in the introduction, so we may expect some differences in the argumentation for $k=3$. Nonetheless, a substantial part of the proof of Theorem~\ref{T6} can be taken over.

{\bf Proof of Theorem~\ref{T1}.}
Let $Q$ be a DFS tree of a $3$-connected graph $G$ nonisomorphic to $K_4$. We will prove that $Q$ contains at least one $3$-contractible edge,
and that $Q$ contains at least two $3$-contractible edges unless $G$ is a prism or a prism plus a single edge (and $Q$ has a special shape).
Set $\mathfrak{S}:=\{V(e):\, e \in E(Q)\}$ and $\mathfrak{R}:=\{V(e):\, e \in E(Q)$, $|V(e) \cap A|=1$ for some $\mathfrak{S}$-end $A\}$.
Let us adopt once more the terminology of small, big, good, bad, and very good fragments introduced in the proof of Lemma~\ref{L2}.

By Lemma~\ref{L2} (which is, in contrast to Theorem~\ref{T7}, true for $k=3$, too), we may assume that there exists a small, bad $T_B$-$\mathfrak{R}$-end $B=\{b\}$.
There exists an edge $e \in E(Q)$ and a $T_A$-end $A$ such that $V(e) \subseteq T_B$ and $|V(e) \cap A|=1$, and as in Claim~1 in the proof of Theorem~\ref{T6}
we see that $|A|=1$, say, $A=\{a\}$. Let $x$ be the vertex in $V(e) \setminus \{a\}$ and let $t$ be the vertex in $T_B \setminus V(e)$.
 As $bt$ is $3$-contractible (because $N_G(b) \setminus \{t\}$ induces a triangle in $G$), we know that $bt \not \in E(Q)$, as otherwise $B$ would not be bad. Hence, as $ax \in E(Q)$,
exactly one of $ba,bx$ is in $E(Q)$. Observe that $a$ is not adjacent to $t$, because $a$ has degree $3$ and a neighbor in $\overline{B}$.

Suppose that $ba \in E(Q)$. Then there exists a $T \in \mathfrak{T}(G)$ containing $b,a$, and a vertex $s$ from $\overline{B}$. 
Consider any $T$-fragment $F$. At least one of $F \cap \overline{B},\overline{F} \cap \overline{B}$ must be empty,
for otherwise one of these sets is an $\{a,s,x\}$-fragment and the other one is an $\{a,s,t\}$-fragment; this is not possible, because $a$ has only one neighbor in $\overline{B}$.
So $\{x\}$ or $\{t\}$ is a $T$-fragment, but the latter is not because $t$ is not adjacent to $a \in T$.
Consequently, $\{x\}$ is a $T$-fragment, and, in fact, an $\mathfrak{S}$-atom. 

Therefore, after possibly having exchanged the names of $a,x$ (and resetting $A$ accordingly), we may assume without loss of generality
that $bx \in E(Q)$; $T_A$ consists of $b,x$, and a vertex $s$ from $\overline{B}$, and we see immediately that $A=\{a\}$ is a small $\mathfrak{R}$-end, too.
The edge $as$ is 3-contractible, because $N_G(a) \setminus \{s\}=\{b,x\}$ induces a complete graph in $G$.

Now if the small $\mathfrak{R}$-end $A$ was bad, $ax$ would be the only edge incident with $a$ in $Q$, since $as$ is $3$-contractible.
Then one of $a$ and $b$ must be the root vertex $r$ of the DFS tree $Q$,
because otherwise $a,b$ would belong to different $x$-branches but are adjacent. If $r=a$, then $b$ and $t$ belong to different $x$-branches but are adjacent; if otherwise $r=b$, then $a$ and $s$ belong to different $x$-branches but are adjacent.
Thus, $A$ is good.

Since $ab \notin E(Q)$ and $ax$ is not $3$-contractible, 
it follows that $Q$ contains the $3$-contractible edge $as$. This gives the first claim; for the second, let us assume that $as$ is the only $3$-contractible edge from $Q$.
We have to prove that $G$ is either the prism or the prism plus a single new edge; to this end, we proceed as in the proof of Theorem~\ref{T7},
and consider an arbitrary $T_{B'}$-$\mathfrak{R}$-end $B' \subseteq \overline{A}$. 
If it was good, we would find a $3$-contractible edge having an
endvertex in common with $B'$, and, thus, distinct from $as$, which gives a contradiction. Thus, $B'$ is bad. If $B'$ was big, it would contain
a very good end in its neighborhood; just as in the proof of Theorem~\ref{T7}, we cannot locate the root of $Q$ properly, which gives a contradiction.

Therefore, $B'$ is a bad, small $\mathfrak{R}$-end, and we may apply to $B'$ the same line of arguments that we applied before to $B$: $T_{B'}$ consists of three vertices $a',x',t'$,
where $A':=\{a'\}$ is a good, small $T_{A'}$-$\mathfrak{R}$-end, and $T_{A'}$ consists of $b',x'$ and a vertex $s'$ from $\overline{B'}$,
where, moreover, $x'a'$ and $x'b'$ are from $E(Q)$ and not $3$-contractible, and $a's'$ is from $E(Q)$ and $3$-contractible.

Since $as$ is the only $3$-contractible edge, we see that $a's'=as$, which implies $a'=s$ and $s'=a$.
Again the situation is symmetric in $a,b,x,s,t$ and $a',b',x',s',t'$, and we may proceed almost literally as in the proof of Theorem~\ref{T7}
by showing first that $x \not= x'$, $r \notin \{x,x'\}$, and $R:=Q[\{a,b,x,a',b',x'\}]$ is a subtree of $Q$, and thus, more precisely, a path $bxaa'x'b'$.
By symmetry, we may assume that $r \in \{a',b'\}$, and, as $b$ cannot be on the $r,t$-path in $Q$, that $t$ must be on the $r,b$-path in $Q$ and hence in $R \cap \overline{A}$.
Different from the more general argument,
we now have two options for $t$: $t$ is either $b'$, or it is $x'$. We consider the corresponding cases separately:

{\bf Case 1.} $t=b'$

Then $t' = b$. The set $X:=\overline{A} \cap \overline{B} \cap \overline{A'} \cap \overline{B'}$ has all its neighbors
in $T_A \cup T_B \cup T_{A'} \cup T_{B'} \setminus (A \cup B \cup A' \cup B')=\{x,x'\}$. Therefore, $X$ is empty,
and so $\overline{A} \cap \overline{B}=\{x'\}$, which implies that $G$ is the prism.

{\bf Case 2.} $t=x'$. 

If there are no further vertices but $a,b,a',b',x,x'=t,t'$ and they are all distinct, then $N_G(t') = \{x,t,b'\}$.
As certified by the $\mathfrak{S}$-end $\{b'\}$ and the edge $b't$, $\{t'\}$ is an $\mathfrak{R}$-end contained in $\overline{A}$.
As we have seen above (for the arbitrarily chosen $B'$), the neighborhood of such a fragment necessarily contains $a'$, contradiction.

If there are further vertices other than those listed in the previous paragraph, they are all from
$X:=\overline{A} \cap \overline{B} \cap \overline{A'} \cap \overline{B'}$. Then $X$ has all its neighbors
in $T_A \cup T_B \cup T_{A'} \cup T_{B'} \setminus (A \cup B \cup A' \cup B')=\{x,x',t'\}=:T$. 
Therefore, $t' \not= x$, and $X$ is, indeed, an $T$-fragment, and $\overline{X}=\{a,b,a',b'\}$.
Since $t'$ has only one neighbor in $\overline{X}$, $F:=X \cup \{t'\}$ is an $\{x,x',b'\}$-fragment (with $\overline{F}=\{a,b,a'\}$).
The $\mathfrak{S}$-end $\{b'\}$ together with the edge $b't$ certifies that $F$ is an $\mathfrak{R}$-fragment,
and, thus, contains an $\mathfrak{R}$-end $B''$, which is in $\overline{A}$. However, the neighborhood of $B''$ does not contain $a'$,
as it should, by what we have proven about the arbitrarily chosen $B'$ above, which gives a contradiction.

Therefore, there are no further vertices but $a,b,a',b',x,x'=t,t'$, and the vertices listed are not all distinct.
This implies $x=t'$. The neighborhoods of $a,b,a',b'$ are determined, so that $E(G)$ is determined up to a possible edge connecting $x$ and $x'$.
If $x$ and $x'$ are not adjacent, we get the prism, and otherwise we get the prism plus a single edge.

This proves the theorem.
\hspace*{\fill}$\Box$

\bigskip

{\bf Address of the authors.}

Institut f\"ur Mathematik \\
Technische Universit\"at Ilmenau \\
Weimarer Stra{\ss}e 25 \\
D--98693 Ilmenau \\
Germany

\end{document}